\documentclass{amsart}
\usepackage{amssymb}

\newtheorem{theorem}{Theorem}

\newtheorem{lemma}[theorem]{Lemma}

\numberwithin{equation}{section}
\input{tcilatex}

\begin{document}
\title[Spherical CR manifolds]{Connected sum of spherical CR manifolds with
positive CR Yamabe constant}
\author{Jih-Hsin Cheng}
\address{Institute of Mathematics, Academia Sinica and National Center for
Theoretical Sciences, Taipei, Taiwan, R.O.C.}
\email{cheng@math.sinica.edu.tw}
\author{Hung-Lin Chiu}
\address{Department of Mathematics, National Tsing-Hua University, Hsinchu,
Taiwan, R.O.C.}
\email{hlchiu@math.nthu.edu.tw}
\subjclass{32V05, 32V20}
\keywords{connected sum, Spherical CR Manifolds, CR Yamabe constant}
\thanks{}

\begin{abstract}
Suppose $M_{1}$ and $M_{2}$ are two closed (compact with no boundary)
spherical CR manifolds with positive CR Yamabe constant. In this note, we
show that the connected sum of $M_{1}$ and $M_{2}$ also admits a spherical
CR structure with positive CR Yamabe constant.
\end{abstract}

\maketitle

\section{\textbf{Introduction and statement of the results}}

In Riemannian geometry, we have the following fact about positive scalar
curvature. Namely, the connected sum of two closed (compact with no
boundary) manifolds of positive scalar curvature has a metric of positive
scalar curvature (see Corollary 3 in Schoen-Yau's paper \cite{SY} or
Gromov-Lawson's paper \cite{GL}). This fact has been generalized to
surgeries in codimension $\geq $ $3$ and discussed in (spin or not spin)
cobordism theory in dimension $\geq $ $5$ (see \cite{GL}). As an interesting
result in conformal geomery, the connected sum of two closed conformally
flat manifolds with positive Yamabe constant is still a conformally flat
manifolds with positive Yamabe constant (see Corollary 5 in \cite{SY}). In
this note, we are going to prove an analogue in CR geometry.

For basic material in $CR$ and pseudohermitian geometry, we refer the reader
to \cite{Lee} or \cite{Web}. Let $(M,J)$ be a closed, strictly pseudoconvex $%
CR$ manifold of dimension 2n+1. Take a contact form $\theta ,$ so we can
talk about $L^{p}$ norm $||\cdot ||_{p},$ Levi metric $|\cdot |,$
subgradient $\nabla _{b}$ and Tanaka-Webster scalar curvature $R$ or $%
R_{J,\theta }$ on the pseudohermitian manifold $(M,J,\theta ).$ Take the
volume form $dV$ $:=$ $\theta \wedge (d\theta )^{n}.$ Then we can write the $%
CR$ Yamabe constant (or invariant) $\mathcal{\lambda }(M,J)$ or $\mathcal{%
\lambda }(M)$ (if $J$ is clear in the context) as%
\begin{equation*}
\mathcal{\lambda }(M,J)=\inf_{u>0}\frac{E_{\theta }(u)}{||u||_{2+2/n}^{2}}
\end{equation*}

\noindent where%
\begin{equation*}
E_{\theta }(u):=\int_{M}[(2+\frac{2}{n})|\nabla _{b}u|^{2}+Ru^{2}]dV
\end{equation*}

\noindent and%
\begin{equation*}
||u||_{2+2/n}^{2}:=(\int_{M}|u|^{2+\frac{2}{n}}dV)^{\frac{n}{n+1}}
\end{equation*}

\noindent (see \cite{JL} for more details).

\bigskip

\textbf{Theorem A}. \textit{Suppose }$(M_{1},J_{1})$\textit{\ and }$%
(M_{2},J_{2})$\textit{\ are two closed, spherical }$CR$\textit{\ manifolds
of dimension }$2n+1$\textit{\ with }$\lambda (M_{k},J_{k})$\textit{\ }$>$%
\textit{\ }$0$\textit{\ for }$k=1,2.$\textit{\ Then their connected sum }$%
M_{1}\#M_{2}$\textit{\ admits a spherical }$CR$\textit{\ structure }$\tilde{J%
}$\textit{\ with }$\lambda (M_{1}\#M_{2},\tilde{J})$\textit{\ }$>$\textit{\ }%
$0.$

\bigskip

The idea of the proof was motivated by the work of O. Kobayashi \cite{Kob}
(note that we do not follow the approach of either Schoen-Yau or
Gromov-Lawson). Kobayashi's short proof is specially suitable for
conformally flat manifolds with positive Yamabe constant. Due to different
nature of cylinders structure between conformal and CR geometries, we modify
the original idea of Kobayashi for the CR case (after we posted our paper on
the ArXiv, we learned from Yun Shi that a similar argument in \cite{SW}
works for the original approach of Kobayashi in the CR case. This fills up
the gap of the proof in \cite{Wang}).

Theorem A is used to construct many examples in the study of positive mass
theorem for 5 dimensional closed, strictly pseudoconvex CR manifolds $M$ (%
\cite{CC}). In \cite{CC}, we assume further $M$ is spin, spherical with
positive CR Yamabe constant. Then we have positive mass theorem for $M.$
According to Theorem A, we have the following examples:

\begin{equation*}
m_{1}(S^{5}/Z_{p_{1}})\#l_{1}(S^{4}\times
S_{(a_{1})}^{1})\#m_{2}(S^{5}/Z_{p_{2}})\#l_{2}(S^{4}\times
S_{(a_{2})}^{1})\#...
\end{equation*}

\bigskip

\noindent $($connected sum of finite number of manifolds such as $%
S^{5}/Z_{p} $ or $S^{4}\times S_{(a_{j})}^{1},$ $a_{j}$ $>$ $1)$ for $p_{j}$
odd (noting that $S^{5}/Z_{2}$ is not spin, but still spherical with
positive CR Yamabe constant). See the end of Section 2 for more details.

\bigskip

\textbf{Acknowledgements. }J.-H. C. would like to thank Jack Lee, Matt
Gursky and Paul Yang for useful discussions and informing him of papers of
Gromov-Lawson and Schoen-Yau. We would also like to thank Kengo Hirachi to
inform us of the related paper \cite{Wang} after our paper was posted on the
ArXiv. J.-H. C. (H.-L.C., resp.) would like to thank the Ministry of Science
and Technology of Taiwan for the grant \# 106- 2115- M- 001- 013. (MOST
106-2115-M-007-017-MY3, resp.).

\bigskip

\section{\textbf{Proof of Theorem A}}

%\subsection{The Heisenberg cylinder}

We first discuss the structure of the Heisenberg cylinder. Let $H_{n}$
denote the Heisenberg group. On $H_{n}\setminus \{0\}$, the dilations $\tau
_{a}(z,t)=(az,a^{2}t),$ $z$ $=$ $(z^{1},$ $...,$ $z^{n})$ $\in $ $C^{n},$ $t$
$\in $ $R,$ for $a>0$ and the CR inversion map $(z,t)\rightarrow (z^{\ast
},t^{\ast })$ defined by%
\begin{equation*}
z^{\ast }=\frac{z}{w},\ \ t^{\ast }=-\frac{t}{|w|^{2}},\ \ \text{where}\
w:=t+i|z|^{2},
\end{equation*}%
are all CR transformations. The standard contact form $\Theta $ on $H_{n}$
reads%
\begin{equation*}
\Theta :=dt+i\sum_{\alpha =1}^{n}(z^{\alpha }dz^{\bar{\alpha}}-z^{\bar{\alpha%
}}dz^{\alpha }).
\end{equation*}%
\noindent Instead of $\Theta $, we consider the contact form $\frac{\Theta }{%
\rho ^{2}}$ where $\rho $ $:=$ $|w|^{1/2}$ $=$ ($|z|^{4}+t^{2})^{1/4}.$ Then
all these maps also preserve the new contact form $\frac{\Theta }{\rho ^{2}}$%
. The space $H_{n}\setminus \{0\}$ together with the new contact form $\frac{%
\Theta }{\rho ^{2}}$ is called \textbf{the Heisenberg cylinder}.
Topologically, $H_{n}\setminus \{0\}$ $=$ $(0,\infty )\times S^{2n}(1)$
where $S^{2n}(1)$ $:=$ $\{\rho =1\}$ $\subset $ $H_{n}.$ For fixed $a>1$,
each slice $[a^{m-1},a^{m}]\times S^{2n}(1)$ is isomorphic to one another as
pseudohermitian manifolds. Consider the quotient space%
\begin{eqnarray}
&&S^{2n}\times S^{1}\text{ or }S^{2n}\times S_{(a)}^{1}\text{ (to indicate
the dependence on }a)  \label{HC} \\
&=&H_{n}\setminus \{0\}/\{\cdots ,\tau _{a^{-1}},1,\tau _{a}.\tau
_{a^{2}},\cdots \}.  \notag
\end{eqnarray}

We want to apply the following results to $S^{2n}\times S_{(a)}^{1}.$ The
proof is similar as for the analogous statements in the Riemannian case
(see, e.g., \cite{LP}). For completeness, we give a proof here.

\bigskip

\begin{lemma}
\label{lem0} Let $(M,J,\theta )$ be a closed pseudohermitian manifold of
dimension $2n+1.$

(1) Suppose the Tanaka-Webster scalar curvature $R_{J,\theta }\geq 0,$ $>0$
somewhere. Then $\mathcal{\lambda }(M,J)>0.$

(2) Suppose $\mathcal{\lambda }(M,J)>0.$ Then there exists a contact form $%
\tilde{\theta}$ such that $R_{J,\tilde{\theta}}>0.$
\end{lemma}

\bigskip

\proof For (1), suppose $\mathcal{\lambda }(M,J)$ $\leq $ $0.$ We can solve
the $CR$ Yamabe equation by a theorem in \cite{JL} to find $u$ $>$ $0$ such
that%
\begin{equation}
\mathcal{-(}2\mathcal{+}\frac{2}{n})\Delta _{b}u+R_{J,\theta }u=\mathcal{%
\lambda }(M,J)u^{\frac{n+2}{n}}.  \label{Eq}
\end{equation}

\noindent Multiplying (\ref{Eq}) by $u$ and integrating give%
\begin{eqnarray*}
0 &\leq &\int_{M}[\mathcal{(}2\mathcal{+}\frac{2}{n})|\nabla
_{b}u|^{2}+R_{J,\theta }u^{2}]dV \\
&=&\mathcal{\lambda }(M,J)\int_{M}u^{2\mathcal{+}\frac{2}{n}}dV\leq 0
\end{eqnarray*}

\noindent since $R_{J,\theta }\geq 0$ and $\mathcal{\lambda }(M,J)\leq 0.$
It follows that $R_{J,\theta }$ $\equiv $ $0.$ Contradicts to $R_{J,\theta }$
$>$ $0$ somewhere. So we conclude $\mathcal{\lambda }(M,J)>0.$

For (2), let $2\leq s$ $<$ $p$ :$=$ $2\mathcal{+}\frac{2}{n},$ critical
exponent. Set $\lambda _{s}$ $:=$ $\inf \{E_{\theta }(\varphi )/||\varphi
||_{s}^{2}:$ $\varphi $ $\in $ $C^{\infty }(M)\}.$ There exists a smooth,
positive solution $\varphi _{s}$ to the subcritical equation%
\begin{equation}
\mathcal{-}p\Delta _{b}\varphi _{s}+R_{J,\theta }\varphi _{s}=\lambda
_{s}\varphi _{s}^{s-1}  \label{2-1}
\end{equation}

\noindent (Folland-Stein space $S_{1}^{2}\subset L^{s}$ is compact). On the
other hand, let $\tilde{\theta}$ $:=$ $\varphi _{s}^{2/n}\theta $. Then by
the transformation law, we have%
\begin{equation}
\mathcal{-}p\Delta _{b}\varphi _{s}+R_{J,\theta }\varphi _{s}=R_{J,\tilde{%
\theta}}\varphi _{s}^{p-1}  \label{2-2}
\end{equation}

\noindent From (\ref{2-1}) and (\ref{2-2}), we have%
\begin{equation}
R_{J,\tilde{\theta}}=\lambda _{s}\varphi _{s}^{s-p}.  \label{2-3}
\end{equation}

\noindent Observe that $\lambda _{s}$ is continuous in $s$ from the left
(cf. Lemma 4.3 in \cite{LP}). So $\lambda _{s}$ $>$ $0$ for $s$ close to $p$
where $\lambda _{p}=\mathcal{\lambda }(M,J)>0$ by assumption. Thus $R_{J,%
\tilde{\theta}}$ $>$ $0$ in view of (\ref{2-3}) for $s$ close to $p.$

\endproof

\bigskip

We compute the Tanaka-Webster scalar curvature $R=R_{\frac{\Theta }{\rho ^{2}%
}}$ on $S^{2n}\times S_{(a)}^{1}$ as follows:%
\begin{equation*}
R_{\frac{\Theta }{\rho ^{2}}}=\frac{n(n+1)|z|^{2}}{2\rho ^{2}}.
\end{equation*}

\noindent Observe that $R_{\frac{\Theta }{\rho ^{2}}}\geq 0$ and $>$ $0$ if $%
z\neq 0.$ From Lemma \ref{lem0} (1) and (2), we can find a contact form $%
\hat{\theta}$ such that the Webster curvature $R_{\hat{\theta}}$, with
respect to this contact form $\hat{\theta}$, is positive on $S^{2n}\times
S^{1}$. Therefore if we consider the lifting $\hat{\Theta}$ of $\hat{\theta}$
by the covering map $H_{n}\setminus \{0\}\rightarrow (S^{2n}\times S^{1},%
\hat{\theta})$ then the Webster curvature $R_{\hat{\Theta}}>c$ for some
positive constant $c$. In addition, $\tau _{a}$ also defines a symmetry on $%
(H_{n}\setminus \{0\},\hat{\Theta})$, which implies that each slice $%
[a^{m-1},a^{m}]\times S^{2n}(1)$ with respect to this new contact form $\hat{%
\Theta}$, instead of $\frac{\Theta }{\rho ^{2}}$, is also isomorphic to one
another as pseudohermitian manifolds.

\bigskip

\begin{lemma}
\label{lem1} \textit{We have }%
\begin{equation}
\lambda ((M_{1},J_{1})\amalg (M_{2},J_{2}))=\min \{\lambda
(M_{1},J_{1}),\lambda (M_{2},J_{2})\},
\end{equation}%
\noindent \textit{provided that both }$\lambda (M_{1},J_{1})$\textit{\ and }$%
\lambda (M_{2},J_{2})$\textit{\ are positive.}
\end{lemma}

\bigskip

\proof Let $M=(M_{1},J_{1})\amalg (M_{2},J_{2})$ and $f=f_{1}\amalg f_{2}$
be a $C^{\infty }$ function on $M$, where $f_{1}$ and $f_{2}$ are $C^{\infty
}$ functions on $M_{1}$ and $M_{2},$ resp.. If we take $f_{2}$ to be zero,
it is easy to see that $\lambda (M)\leq \lambda (M_{1},J_{1})$. Similarly,
taking $f_{1}$ to be zero, we have $\lambda (M)\leq \lambda (M_{2},J_{2})$.

On the other hand, suppose we choose $f=f_{1}\amalg f_{2}$ such that 
\begin{equation}
\int_{M}|f|^{2+\frac{2}{n}}dV=1.  \label{norcon}
\end{equation}%
We compute 
\begin{equation*}
\begin{split}
& (2+\frac{2}{n})\int_{M}|\nabla _{b}f|^{2}dV+\int_{M}Rf^{2}dV \\
=& \sum_{j=1}^{2}(2+\frac{2}{n})\int_{M_{j}}|\nabla
_{b}f_{j}|^{2}dV_{j}+\int_{M_{j}}R_{j}f_{j}^{2}dV_{j} \\
\geq & \sum_{j=1}^{2}\lambda (M_{j},J_{j})\left( \int_{M_{j}}|f_{j}|^{2+%
\frac{2}{n}}dV_{j}\right) ^{\frac{n}{n+1}} \\
\geq & \sum_{j=1}^{2}\lambda (M_{j},J_{j})\left( \int_{M_{j}}|f_{j}|^{2+%
\frac{2}{n}}dV_{j}\right) =\sum_{j=1}^{2}\lambda (M_{j},J_{j})\alpha _{j},
\end{split}%
\end{equation*}%
where $\alpha _{j}=\int_{M_{j}}|f_{j}|^{2+\frac{2}{n}}dV_{j}>0$, and by (\ref%
{norcon}), $\alpha _{1}+\alpha _{2}=1$. This shows that $\lambda (M)\geq
\min \{\lambda (M_{1},J_{1}),\lambda (M_{2},J_{2})\}$.

\endproof

\bigskip

\proof\textbf{(of Theorem A)} Let $M_{j},j=1,2$, be two differentiable
manifolds and let $M=M_{1}\amalg M_{2}$ be the disjoint union of $M_{1}$ and 
$M_{2}$. Fix $p_{j}$ $\in $ $M_{j},$ $j=1,2.$ We take off two small balls
around $p_{1}$ and $p_{2}$, and then attach a cylinder, which is
topologically the product of a line segment and $S^{2n}$. The new manifold
obtained in this way is called the connected sum of $M_{1}$ and $M_{2}$
denoted by $M_{1}\#M_{2}$. If, in addition, assume that $M_{j},j=1,2$, are
two spherical CR manifolds with pseudohermitian structures $(J_{j},\theta
_{j})$, then we will use (part of) the Heisenberg cylinder to glue them
together in order that the result manifold is also spherical.

We can find contact forms $\tilde{\theta}_{j}$ on $M_{j}\setminus \{p_{j}\},$
$j=1,2,$\ such that it is part of the Heisenberg cylinder on a punched
neighborhood of $p_{j}$. Precisely, we can choose $\tilde{\theta}_{j}$ such
that 
\begin{equation*}
M_{j}\setminus \{p_{j}\}=\tilde{M}_{j}\cup \lbrack 1,\infty )\times S^{2n}(1)
\end{equation*}%
or equivalently 
\begin{equation*}
M_{j}\setminus \{p_{j}\}=\tilde{M}_{j}\cup (0,1]\times S^{2n}(1)
\end{equation*}

For convenience, we write 
\begin{equation*}
\left( M\setminus \{p_{1},p_{2}\},\tilde{J},\tilde{\theta}\right) =\Big(%
\lbrack 1,\infty )\times S^{2n}(1)\cup \tilde{M}_{1}\Big)\amalg \Big(\tilde{M%
}_{2}\cup (0,1]\times S^{2n}(1)\Big),
\end{equation*}%
where 
\begin{equation*}
\tilde{J}|_{M_{j}\setminus \{p_{j}\}}=J_{j},\ \ \ \tilde{\theta}%
|_{M_{j}\setminus \{p_{j}\}}=\tilde{\theta}_{j}.
\end{equation*}%
Fix $a>1$ and $l\in N$, a positive integer. Let $\bar{M}=M_{1}\#M_{2}$ be
the connected sum of $M_{1}$ and $M_{2}$ obtained by cutting ends $%
(a^{l},\infty )\times S^{2n}(1)$ and $(0,a^{-l})\times S^{2n}(1)$ out of $%
M\setminus \{p_{1},p_{2}\}$ and gluing the left two parts by means of the
dilation $\tau _{a^{-l}}:[1,a^{l}]\times S^{2n}(1)\rightarrow \lbrack
a^{-l},1]\times S^{2n}(1)$. We hence get a new spherical CR manifold $(\bar{M%
},J_{l},\theta _{l})$ such that 
\begin{equation}
\bar{M}=\tilde{M}\cup \lbrack 1,a^{l}]\times S^{2n}(1),\ \ \text{where}\ 
\tilde{M}=\tilde{M}_{1}\cup \tilde{M}_{2};
\end{equation}%
\begin{equation*}
J_{l}|_{\tilde{M}_{j}}=\tilde{J}=J_{j},\ \ \theta _{l}|_{\tilde{M}_{j}}=%
\tilde{\theta}_{j}\ \ \ \text{and}\ \theta _{l}|_{[1,a^{l}]\times S^{2n}(1)}=%
\frac{\Theta }{\rho ^{2}}.
\end{equation*}%
Moreover, instead of $\frac{\Theta }{\rho ^{2}}$, we can choose $\theta _{l}$
such that $\theta _{l}|_{[1,a^{l}]\times S^{2n}(1)}=\hat{\Theta}$ so that
the Webster curvature is strictly positive on each part $[a^{m-1},a^{m}]%
\times S^{2n}(1)$, for $m=1\cdots l$.

Recall that 
\begin{equation*}
\lambda (\bar{M},J_{l})=\inf_{f>0}\frac{(2+\frac{2}{n})\int_{\bar{M}}|\nabla
_{b}f|^{2}dV_{l}+\int_{\bar{M}}R_{l}f^{2}dV_{l}}{\left( \int_{\bar{M}}|f|^{2+%
\frac{2}{n}}dV_{l}\right) ^{\frac{n}{n+1}}},
\end{equation*}%
where $dV_{l}=\theta _{l}\wedge (d\theta _{l})^{n}$. So, take a positive
function $f_{l}\in C^{\infty }(\bar{M})$ such that 
\begin{equation}
(2+\frac{2}{n})\int_{\bar{M}}|\nabla _{b}f_{l}|^{2}dV_{l}+\int_{\bar{M}%
}R_{l}f_{l}^{2}dV_{l}\leq \lambda (\bar{M},J_{l})+\frac{1}{l}  \label{1}
\end{equation}%
and 
\begin{equation}
\int_{\bar{M}}|f_{l}|^{2+\frac{2}{n}}dV_{l}=1.  \label{2}
\end{equation}

\bigskip

\begin{lemma}
\label{lem2} \textit{There is an integer }$m\in \{1,\cdots ,l\}$\textit{\
such that }%
\begin{equation}
\int_{\lbrack a^{m-1},a^{m}]\times S^{2n}(1)}(|\nabla
_{b}f_{l}|^{2}+f_{l}^{2})dV_{l}\leq \frac{A}{l},  \label{3}
\end{equation}%
\textit{where }$A$\textit{\ is a constant independent of }$l\in N$\textit{.}
\end{lemma}

\bigskip

\proof From (\ref{1}), we have 
\begin{equation}
\begin{split}
& (2+\frac{2}{n})\int_{[1,a^{l}]\times S^{2n}(1)}|\nabla
_{b}f_{l}|^{2}dV_{l}+\int_{[1,a^{l}]\times S^{2n}(1)}R_{l}f_{l}^{2}dV_{l} \\
& \leq \lambda (\bar{M},J_{l})+\frac{1}{l}+\int_{\tilde{M}%
}(-R_{l}f_{l}^{2})dV_{l}.
\end{split}
\label{4}
\end{equation}%
On the other hand, using (\ref{2}) and H\"{o}lder inequality, we have 
\begin{equation}
\begin{split}
\int_{\tilde{M}}-R_{l}f_{l}^{2}dV_{l}& \leq \left( \int_{\tilde{M}%
}|R_{l}|^{n+1}dV_{l}\right) ^{\frac{1}{n+1}}\left( \int_{\tilde{M}}|f_{l}|^{%
\frac{2(n+1)}{n}}dV_{l}\right) ^{\frac{n}{n+1}} \\
& \leq \left( \max_{\tilde{M}}|\tilde{R}|\right) \left( vol(\tilde{M}%
)\right) ^{\frac{1}{n+1}},
\end{split}
\label{5}
\end{equation}%
where $\tilde{R}$ is the Webster curvature with respect to $(\tilde{J},%
\tilde{\theta})$. Substituting (\ref{5}) into (\ref{4}) and noticing that $%
R_{l}$ has an uniform lower bound $c>0$ on the Heisenberg cylinder, we have 
\begin{equation}
\begin{split}
& (2+\frac{2}{n})\int_{[1,a^{l}]\times S^{2n}(1)}|\nabla
_{b}f_{l}|^{2}dV_{l}+c\int_{[1,a^{l}]\times S^{2n}(1)}f_{l}^{2}dV_{l} \\
& \leq \lambda (S^{2n+1})+\frac{1}{l}+A_{1},
\end{split}
\label{6}
\end{equation}%
where $\lambda (S^{2n+1})$ is the Yamabe constant of the standard sphere and 
\begin{equation*}
A_{1}=\left( \max_{\tilde{M}}|\tilde{R}|\right) \left( vol(\tilde{M})\right)
^{\frac{1}{n+1}},
\end{equation*}%
which is independent of $l$. Let $C=\min \{2+\frac{2}{n},c\}$ and $A=\frac{%
\lambda (S^{2n+1})+1+A_{1}}{C}$, and let $m\in \{1,\cdots ,l\}$ be chosen so
that the energy $\int_{[a^{m-1},a^{m}]\times S^{2n}(1)}(|\nabla
_{b}f_{l}|^{2}+f_{l}^{2})dV_{l}$ on the interval $[a^{m-1},a^{m}]$ is the
smallest among those on all the $l$ intervals. Then this energy satisfies
the assertion of this lemma. We have completed the proof.

\endproof

\bigskip

Now we cut off $\bar{M}$ on the section $\{a^{\frac{2m-1}{2}}\}\times
S^{2n}(1)$, where $a^{\frac{2m-1}{2}}=\sqrt{a^{m-1}a^{m}}$, and attach
respectively the two cylinders $[a^{\frac{2m-1}{2}},\infty )\times S^{2n}(1)$
and $(0,a^{\frac{2m-1-2l}{2}}]\times S^{2n}(1)$ to it, precisely, to $[1,a^{%
\frac{2m-1}{2}}]\times S^{2n}(1)$ and $[a^{\frac{2m-1-2l}{2}},1]\times
S^{2n}(1)$. Then we obtain again the manifold 
\begin{equation*}
\begin{split}
& (M\setminus \{p_{1},p_{2}\},\tilde{J},\tilde{\theta}) \\
=& \left( \bar{M}\setminus (\{a^{\frac{2m-1}{2}}\}\times S^{2n}(1))\right)
\cup \lbrack a^{\frac{2m-1}{2}},\infty )\times S^{2n}(1)\cup (0,a^{\frac{%
2m-1-2l}{2}}]\times S^{2n}(1).
\end{split}%
\end{equation*}%
We think of the function $f_{l}$ as defined on $\bar{M}\setminus (\{a^{\frac{%
2m-1}{2}}\}\times S^{2n}(1))$, and extend it to the whole space $(M\setminus
\{p_{1},p_{2}\},\tilde{J},\tilde{\theta})$ as follows: Let $F_{l}$ be the $%
C^{\infty }$ function on $M\setminus \{p_{1},p_{2}\}$ such that 
\begin{equation}
F_{l}=\left\{ 
\begin{array}{ll}
f_{l} & \text{on}\ \bar{M}\setminus \{a^{\frac{2m-1}{2}}\}\times S^{2n}(1)
\\ 
\chi _{1}f_{l} & \ \text{on}\ [a^{\frac{2m-1}{2}},\infty )\times S^{2n}(1)
\\ 
\chi _{2}f_{l} & \ \text{on}\ (0,a^{\frac{2m-1-2l}{2}}]\times S^{2n}(1),%
\end{array}%
\right.
\end{equation}%
where both $\chi _{1}$ and $\chi _{2}$ are cut off functions on the cylinder
defined by 
\begin{equation*}
\chi _{1}(t,x)=\left\{ 
\begin{array}{ll}
1, & t\leq a^{\frac{2m-1}{2}}; \\ 
0, & t\geq a^{m},%
\end{array}%
\right.
\end{equation*}%
and 
\begin{equation*}
\chi _{2}(t,x)=\left\{ 
\begin{array}{ll}
1, & t\geq a^{\frac{2m-1-2l}{2}}; \\ 
0, & t\leq a^{m-1-l}.%
\end{array}%
\right.
\end{equation*}%
Notice that each part $[a^{m-1},a^{m}]\times S^{2n}(1)$ is isomorphic to $%
S^{2n}\times S^{1}$, so that $\chi _{1}$ and $\chi _{2}$ can be chosen so
that $|\nabla _{b}\chi _{1}|$ and $|\nabla _{b}\chi _{2}|$ have upper bound $%
c_{1}$ and $c_{2}$ independent of $m$. Now it is easy to see from (\ref{1})
and (\ref{3}) 
\begin{equation*}
(2+\frac{2}{n})\int_{M\setminus \{p_{1},p_{2}\}}|\nabla _{b}F_{l}|^{2}d%
\tilde{V}+\int_{M\setminus \{p_{1},p_{2}\}}\tilde{R}F_{l}^{2}d\tilde{V}\leq
\lambda (\bar{M},J_{l})+\frac{B}{l},
\end{equation*}%
where $d\tilde{V}=\tilde{\theta}\wedge (d\tilde{\theta})^{n}$ and $\tilde{R}$
is the Webster curvature with respect to $(\tilde{J},\tilde{\theta})$, and $%
B $ is a constant independent of $l$. Obviously from (\ref{2}) 
\begin{equation*}
\int_{M\setminus \{p_{1},p_{2}\}}|F_{l}|^{2+\frac{2}{n}}d\tilde{V}>1.
\end{equation*}%
Therefore we have 
\begin{equation}
\inf {\frac{(2+\frac{2}{n})\int_{M\setminus \{p_{1},p_{2}\}}|\nabla
_{b}F|^{2}d\tilde{V}+\int_{M\setminus \{p_{1},p_{2}\}}\tilde{R}F^{2}d\tilde{V%
}}{\left( \int_{M\setminus \{p_{1},p_{2}\}}|F|^{2+\frac{2}{n}}d\tilde{V}%
\right) ^{\frac{n}{n+1}}}}\leq \lambda (\bar{M},J_{l})+\frac{B}{l},
\label{7}
\end{equation}%
where the infimum is taken over all nonnegative $C^{\infty }$ function $F$
with compact support. It follows from the choice of the contact form $\tilde{%
\theta}$ that the left hand side of (\ref{7}) is equal to $\lambda (M)$, in
which $M=M_{1}\amalg M_{2}$. Therefore, by lemma \ref{lem1}, if $l$ is large
enough, we have $\lambda (\bar{M},J_{l})>0$.

\endproof

\bigskip

\textbf{Examples.} Let $p$ be a positive integer. Let $Z_{p}$ denote the
cyclic group generated by the following $(n+1)\times (n+1)$ diagonal matrix:%
\begin{equation*}
\left( 
\begin{array}{cccc}
e^{2\pi i/p} & 0 & .. & 0 \\ 
0 & e^{2\pi i/p} & .. & : \\ 
: & .. & .. & 0 \\ 
0 & .. & 0 & e^{2\pi i/p}%
\end{array}%
\right)
\end{equation*}%
\noindent acting on $C^{n+1}$ by multiplying by $e^{2\pi i/p}$. It is clear
that $Z_{p}$ leaves the unit sphere $S^{2n+1}$ $\subset $ $C^{n+1}$
invariant and preserves the standard $CR$ structure on $S^{2n+1}$. Moreover, 
$Z_{p}$ acts on $S^{2n+1}$ freely. We then have the quotient manifold $%
S^{2n+1}/Z_{p}$ which is closed, spherical and has $\lambda (S^{2n+1}/Z_{p})$
\TEXTsymbol{>} $0$ (since $Z_{p}$ leaves the standard contact form for $%
S^{2n+1}$ invariant, $S^{2n+1}/Z_{p}$ has the same Tanaka-Webster scalar
curvature as $S^{2n+1},$ which is a positive constant$).$ On the other hand,
we learn from the paragraph before Lemma \ref{lem1} that $S^{2n}\times S^{1}$
(see (\ref{HC})) is also a closed, spherical CR manifold with $\lambda
(S^{2n}\times S^{1})$ $>$ $0.$

Now according to Theorem A, we have the following closed, spherical CR
manifolds with positive CR Yamabe constant:%
\begin{equation*}
m_{1}(S^{2n+1}/Z_{p_{1}})\#l_{1}(S^{2n}\times
S_{(a_{1})}^{1})\#m_{2}(S^{2n+1}/Z_{p_{2}})\#l_{2}(S^{2n}\times
S_{(a_{2})}^{1})\#...
\end{equation*}

\noindent (connected sum of finite number of manifolds such as $%
S^{2n+1}/Z_{p_{j}},$ $p_{j}$ being positive integers, or $S^{2n}\times
S_{(a_{j})}^{1},$ $a_{j}$ $>$ $1$) where $m_{1},$ $l_{1},$ $m_{2},$ $l_{2}$
are nonnegative integers. In application to construct 5-dimensional CR
manifolds for positive mass theorem to hold (see \cite{CC}), we need to
restrict $p$ to be odd in order for $S^{5}/Z_{p}$ to be spin.

\end{document}